\documentclass{article}

\usepackage{amssymb}
\usepackage{amsmath}
\usepackage{amsthm}

\newtheorem{theorem}{Theorem} \newtheorem{lemma}{Lemma}[section]
\newtheorem{propo}{Proposition}[section]

\newtheorem{corol}{Corollary}[section] \newtheorem{defin}{Definition}[section]

  \newcommand{\ep}{\varepsilon}
  
 \newcommand{\e}{\ep} 

 \newcommand{\Z}{\mathbb{Z}}
 \newcommand{\C}{\mathbb{C}}

\newcommand{\Mat}{\mbox{Mat}} \newcommand{\Aut}{\mbox{Aut}}
\newcommand{\rk}{\mbox{rk}} \newcommand{\tr}{\mbox{tr}}
\newcommand{\Tr}{\mbox{Tr}}
\newcommand{\cM}{\mathcal{M}}\newcommand{\cD}{\mathcal{D}}
\newcommand{\cN}{\mathcal{N}} \newcommand{\cP}{\mathcal{P}}

\newcommand{\cA}{\mathcal{A}}
\newcommand{\pru}{\Z_{(2)}}  
\newcommand{\dya}{\hat{\Z}_{(2)}}
\newcommand{\haar} {\mu_{\mbox {haar}}} \newcommand{\supp}{\mbox{supp}}
\title{Lamplighter groups and von Neumann`s continuous regular ring\footnote{
\, Research partly sponsored by MTA Renyi ``Lendulet'' Groups
and Graphs Research Group}}

\author{G\'abor Elek}

\begin{document}
\maketitle

\begin{abstract}
Let $\Gamma$ be a discrete group. Following Linnell and Schick one can
define a continuous ring $c(\Gamma)$ associated with $\Gamma$.
They proved that if the Atiyah Conjecture holds for a
 torsion-free group $\Gamma$, 
then $c(\Gamma)$ is
a skew field. Also, if $\Gamma$ has torsion and the Strong Atiyah
 Conjecture holds for $\Gamma$, then
$c(\Gamma)$ is a matrix ring over a skew field. 
The simplest example when the Strong Atiyah Conjecture
fails is the lamplighter group $\Gamma=\Z_2\wr \Z$. It is known that
 $\C(\Z_2\wr \Z)$ does not even have a
classical ring of quotients. Our main result is that if $H$ is amenable, then
$c(\Z_2\wr H)$ is isomorphic to a continuous ring
constructed by John von Neumann in the $1930's$.
\end{abstract}
\tableofcontents
\newpage
\section{Introduction}
Let us consider $\Mat_{k\times k}(\C)$ the algebra of $k$ by $k$ matrices 
over the complex field.
This ring is a unital $*$-algebra with respect to the complex transposes. 
For each element $A\in  \Mat_{k\times k}(\C)$ one can define $A^*$ satisfying
 the following properties.
\begin{itemize}
\item $(\lambda A)^*=\overline{\lambda}A^*$
\item $(A+B)^*=A^*+B^*$
\item $(AB)^*=B^*A^*$
\item ${0}^*=0$, $1^*=1$
\end{itemize}
Also, each element has a normalized rank $\rk(A)=\mbox{Rank} (A)/k$
 with the following properties.
\begin{itemize}
\item $\rk(0)=0, \rk(1)=1, $
\item $\rk(A+B)\leq \rk(A)+\rk(B)$
\item $\rk(AB)\leq \min\{\rk(A),\rk(B)\}$
\item $\rk(A^*)=\rk(A)$
\item If $e$ and $f$ are orthogonal idempotents then $\rk(e+f)=\rk(e)+\rk(f)$.
\end{itemize}
The ring $\Mat_{k\times k}(\C)$ has an algebraic property that von Neumann 
called regularity:
Any principal left-(or right) ideal can be generated by an idempotent. 
Furthermore, among these generating idempotents there is a 
unique projection (that is  $\Mat_{k\times k}(\C)$ is a 
$*$-regular ring). In a von Neumann regular ring any non-zerodivisor is 
necessarily invertible.
One can also observe that
the algebra of matrices is proper, that is $\sum_{i=1}^n a_ia^*_i=0$ 
implies
that all the matrices $a_i$ are zero. One should note that if
$R$ is a $*$-regular ring with a rank function, then the rank
extends to $\Mat_{k\times k}(R)$ \cite{halperin}, where the extended rank has the
same property as $\rk$ except that the rank of the identity is $k$.

\noindent
One can immediately see
 that the rank function defines a metric $d(A,B)\equiv\rk(A-B)$
on any algebra with a rank, and the matrix algebra is complete with
respect to this metric. These complete $*$-regular algebras are called
continuous $*$-algebras
 (see \cite{goodearl} for an extensive study of continuous rings).
Note that for the matrix algebras the possible values of the rank functions 
are $0,1/k,2/k,\dots,1$. John von Neumann observed that there are some
interesting examples of infinite dimensional continuous $*$-algebras, where
the rank function can take any real values in between $0$ and $1$.
His first example was purely algebraic.

\vskip 0.1in
\noindent
{\bf Example 1.} Let us consider the following sequence of diagonal
embeddings.
$$ \C\to \Mat_{2\times 2} (\C)\to  \Mat_{4\times 4} (\C)
 \to \Mat_{8\times 8}  (\C)\to
\dots$$
One can observe that all the embeddings are preserving the rank and the
$*$-operation.
Hence the direct limit $\varinjlim \Mat_{2^k\times 2^k} (\C)$ is a $*$-regular ring
with a proper rank function. The addition, multiplication, the $*$-operation
and the rank function can be extended to the metric completion $\cM$ of
the direct limit ring. The resulting algebra $\cM$ 
is a simple, proper, continuous $*$-algebra,
where the rank function can take all the values on the unit interval.

\vskip 0.1in
\noindent
{\bf Example 2.} Consider a finite, tracial von Neumann algebra $\cN$ with 
trace function $\tr_{\cN}$. Then $\cN$ is a $*$-algebra equipped with a 
rank function.
If $P$ is a projection, then $\rk_{\cN}(P)=\tr_{\cN}(P)$. For a general 
element $A\in\cN$, 
$\rk_{\cN} (A)= 1-\lim_{t\to \infty} \int_0^t \tr_{\cN}(E_\lambda)d\lambda$,
where $\int^\infty_0 E_\lambda\,d\lambda$ is the spectral decomposition of $A^*A$.
In general, $\cN$
is not regular, but it has the Ore property with respect to its zero divisors.
The Ore localization of $\cN$ with respect to its
non-zerodivisors  is called the algebra of
affiliated operators and denoted by $U(\cN)$.
These algebras are also proper continuous $*$-algebras \cite{berberian2}.
The rank of an element $A\in U(\cN)$ is given by the trace of the projection
generating the principal ideal $U(\cN)A$. It is important to note, that
$U(\cN)$ is the rank completion of $\cN$ (Lemma 2.2 (\cite{thom}).

\noindent
Linnell and Schick observed \cite {linnellschick}
that if $X$ is a subset of a proper $*$-regular
algebra  $R$, then there exists a smallest $*$-regular subalgebra containing
$X$, the $*$-regular closure.
 Now let $\Gamma$ be a countable group and $\C\Gamma$ be its complex group 
algebra. Then one can consider the natural embedding of the
group algebra to its group von Neumann algebra $\C\Gamma\to \cN$$ \Gamma$.
Let $U(\Gamma)$ denote the Ore localization of $\cN (\Gamma)$ and the
embedding $\C\Gamma\to U(\Gamma)$. Since $U(\Gamma)$ is a proper $*$-regular
ring, one can consider the smallest $*$-algebra $\cA(\Gamma)$ in $U(\Gamma)$ 
containing
$\C(\Gamma)$. Let $c(\Gamma)$ be the completion of
the algebra $\cA$ above. It is a continuous $*$-algebra \cite{goodearl}.
 Of course, if the rank function has only finitely many
values in $\cA$, then $c(\Gamma)$ equals to $\cA(\Gamma)$.
 Note that if $\C\Gamma$ is
embedded into a continuous $*$-algebra $T$, then one can still define 
$c_T(\Gamma)$ as the smallest continuous ring containing $\C\Gamma$.
 In \cite{elekconnes} we 
proved that if $\Gamma$ is amenable, $c(\Gamma)=c_T(\Gamma)$ for any embedding 
$\C\Gamma \to T$ associated to sofic representations of $\Gamma$, 
hence $c(\Gamma)$ can be viewed as a canonical object.
Linnell and Schick calculated 
the algebra $c(\Gamma)$ for several groups, where the rank function has only 
finitely many values on $\cA$. They proved that the following 
results:
\begin{itemize} 
\item If $\Gamma$ is torsion-free and the Atiyah Conjecture holds
for $\Gamma$, then $c(\Gamma)$ is a skew-field. This is the case, when
$\Gamma$ is amenable and $\C\Gamma$ is a domain. Then $c(\Gamma)$ is the
Ore localization of $\C\Gamma$. If $\Gamma$ is the free group of $k$ generators,
then $c(\Gamma)$ is the Cohen-Amitsur free skew field of $k$ generators.
The  Atiyah Conjecture for a torsion-free group means that the
rank of an element in $\Mat_{k\times k}(\C\Gamma)\subset \Mat_{k\times k}
 (U(\cN(\Gamma)))$
is an integer.
\item If the orders of the finite subgroups of $\Gamma$ are bounded and the
 Strong Atiyah Conjecture holds for $\Gamma$, then $R(\Gamma)$ is
a finite dimensional matrix ring over some skew field. In this case the
Strong Atiyah Conjecture means that the ranks of an element in 
$\Mat_{k\times k}(\C\Gamma)\subset \Mat_{k\times k} (U(\cN(\Gamma)))$ is in 
the abelian
group $\frac{1}{\mbox{lcm}(\Gamma)} \Z$, where $\mbox{lcm}(\Gamma)$ indicates
the least common multiple of the orders of the finite subgroups of $\Gamma$.
\end{itemize}
The lamplighter group $\Gamma=\Z_2\wr \Z$ has finite subgroups of arbitrarily
large orders. 
Also, although $\Gamma$ is amenable, $\C\Gamma$ does 
not satisfy the Ore condition
with respect to its non-zerodivisors \cite{linnell}. In other words, it 
has no classical
ring of quotients. The goal of this paper is to calculate $c(\Z_2\wr \Z)$ and 
even $c(\Z_2\wr H)$, 
where $H$ is a countably infinite amenable group. 
\begin{theorem}\label{t1}
If $H$ is a countably infinite amenable group, then $c(\Z_2\wr H)$ is 
the simple continuous
ring $\cM$ of von Neumann.
\end{theorem}
\section{Crossed Product Algebras}
In this section we recall the notion of
crossed product algebras and the group-measure space construction of
Murray and von 
Neumann. Let $\cA$ be a unital, commutative $*$-algebra and
$\phi:\Gamma\to \Aut(\cA)$ be a representation of the countable group
$\Gamma$ by $*$-automorphisms. The associated crossed product
algebra $\cA\rtimes \Gamma$ is defined the following way.
The elements of $\cA\rtimes \Gamma$ are the finite formal sums
$$\sum_{\gamma\in\Gamma}a_\gamma\cdot\gamma\,,$$
where $a_\gamma \in\cA$.
The multiplicative structure is given by
$$\delta\cdot a_\gamma=\phi(\delta)(a_\gamma)\cdot \delta\,.$$
The $*$-structure is defined by $\gamma^*=\gamma^{-1}$ and
$(\gamma\cdot a)^*=a^*\cdot \gamma^{-1}$.
Note that
$$(\delta\cdot a_\gamma)^*=(\phi(\delta)a_\gamma\cdot \delta)^*=
\delta^*\cdot\phi(\delta)a_\gamma^*=
\phi(\delta^{-1})\phi(\delta) a_\gamma^*\cdot \delta^{-1}=a_\gamma^*\cdot
\delta^*\,.$$

\noindent
Now let $(X,\mu)$ be a probability measure space and $\tau:\Gamma
\curvearrowright X$ be a measure preserving action of a countable group 
$\Gamma$ on $X$. Then we have
a $*$-representation $\hat{\tau}$ of $\Gamma$ in $\Aut (L^\infty(X,\mu))$, 
where $L^\infty(X,\mu)$ is the commutative $*$-algebra of bounded
measurable functions on $X$ (module zero measure perturbations).
$$\hat{\tau}(\gamma)(f)(x)=f(\tau(\gamma^{-1})(x))\,.$$

Let $\mathcal{H}=l^2(\Gamma,L^2(X,\mu))$ be 
the Hilbert-space of 
$L^2(X,\mu)$-valued functions on $\Gamma$.
That is, each element of $\mathcal{H}$ can be written in the form
of $$\sum_{\gamma \in \Gamma} b_\gamma \cdot \gamma\,,$$
where $\sum_{\gamma \in \Gamma} \|b_\gamma\|^2<\infty\,.$
Then we have a representation $L$ of $L^\infty(X,\mu))\rtimes \Gamma$ on
$l^2(\Gamma,L^2(X,\mu))$ by
$$L(\sum_{\gamma\in\Gamma} a_\gamma\cdot\gamma)
(\sum_{\delta\in\Gamma} b_\delta\cdot\delta)=
\sum_{\delta\in\Gamma}(\sum_{\gamma\in\Gamma} 
a_\gamma(\hat{\tau}(\gamma)(\beta_\delta))\cdot \gamma\delta)\,.$$

\noindent
Note that $L(\sum_{\gamma\in\Gamma} a_\gamma\cdot\gamma)$ is always
 a bounded operator.
A trace is given on $L^\infty(X,\mu))\rtimes \Gamma$ by
$$\Tr(S)=\int_X a_1(x) d\mu(x)\,.$$
The weak operator closure
 of $L(L_c^\infty(X,\mu))\rtimes \Gamma)$ in 
$B\left(l^2(\Gamma,L^2(X,\mu))\right)$
 is the von Neumann
algebra $\cN(\tau)$ associated to the action. Here $L^\infty_c(X,\mu)$ denotes
the subspace of functions in $L^\infty(X,\mu)$ having only countable many values.

\noindent
Note that
one can extend $\Tr$ to $\Tr_{\cN(\tau)}$ on the von Neumann algebra to make it
a tracial von Neumann algebra.

\noindent
We will denote by $c(\tau)$ the smallest continuous algebra in
$U(\cN(\tau))$ containing
$L^\infty_c(X,\mu)\rtimes\Gamma$. One should note that the 
weak closure of $L^\infty_c(X,\mu)\rtimes\Gamma$ in
$B\left(l^2(\Gamma,L^2(X,\mu))\right)$ is the same as the weak closure of
$L^\infty(X,\mu)\rtimes\Gamma$. Hence our definition for the von Neumann
algebra of an action coincides with the classical definition. On the
other hand, $c(L^\infty_c(X,\mu)\rtimes\Gamma)$ is smaller
than $c(L^\infty(X,\mu)\rtimes\Gamma)$. 

\section{The Bernoulli Algebra} \label{bernoulli}
Let $H$ be a countable group.
Consider the Bernoulli shift space $B_H:= \prod_{h\in H} \{0,1\}$ with
the usual product measure $\nu_H$. The probability measure preserving action
$\tau_H:H\curvearrowright (B_H,\nu_H)$ is defined by
$$\tau_H(\delta)(x)(h)=x(\delta^{-1} h)\,,$$
where $x\in B_H$, $\delta,h\in H$.
Let $\cA_H$ be the commutative $*$-algebra
of functions that depend only on finitely many
coordinates of the shift space. It is well-known that the
Rademacher functions
$\{R_S\}_{S\subset H,\,|S|<\infty}$ form a basis in $\cA_H$, where
$$R_S(x)=\prod_{\delta\in S} exp(i\pi x(\delta))\,.$$
The Rademacher functions with respect to the pointwise multiplication 
form an Abelian group isomorphic to $\oplus_{h\in H} Z_2$ the
Pontrjagin dual of the compact group $B_H$ satisfying
\begin{itemize}
\item $R_S R_{S'}=R_{S\triangle S'}$
\item $\int_{B_H} R_S \, d\nu=0$\,, if $|S|>0$
\item $R_\emptyset=1$.
\end{itemize}

\noindent
The group $H$ acts on $\cA_H$ by
$$\hat{\tau}_H(\delta)(f)(x)=f\left(\tau_H(\delta^{-1})(x)\right)\,.$$
Hence,
$$\hat{\tau}_H(\delta)R_S=R_{\delta S}\,.$$
Therefore, the elements of $\cA_H\rtimes H$ can be uniquely written as
in the form of the finite sums
$$\sum_{\delta}\sum_S c_{\delta,S} R_S\cdot \delta\,,$$
where $\delta \cdot R_S=R_{\delta S}\cdot \delta\,.$

\noindent
Now let us turn our attention to the group algebra $\C(\Z_2\wr H)$.
For $\delta\in H$, let $t_\delta$ be the generator in $\sum_{h\in H}Z_2$
belonging to the $\delta$-component. Any element of $\C(\Z_2\wr H)$ can be
written in a unique way as a finite sum
$$\sum_\delta\sum_S c_{\delta,S} t_S\cdot \delta\,,$$
where $t_S=\prod_{s\in S} t_s$, $\delta\cdot t_S=t_{\delta S}$, $t_St_{S'}=
t_{S\triangle S'}\,.$
Also note that
$$\Tr(\sum_\delta\sum_S c_{\delta,S} t_S\cdot \delta)=c_{1,\emptyset}\,.$$
Hence we have the following proposition.
\begin{propo}
There exists a trace preserving $*$-isomorphism \\$\kappa:\C(\Z_2\wr H)\to
\cA_H\rtimes H$ such that
$$\kappa(\sum_\delta\sum_S c_{\delta,S} t_S\cdot\delta)=
\sum_\delta\sum_S c_{\delta,S} R_S\cdot\delta\,.$$
\end{propo}

\noindent
Recall that if $A\subset \cN_1$, $B\subset \cN_1$ are weakly dense
*-subalgebras in finite tracial von Neumann algebras $\cN_1$ and $\cN_2$ and
$\kappa:A\to B$ is a trace preserving $*$-homomorphism, then
$\kappa$ extends to a trace preserving isomorphism between the
von Neumann algebras themselves (see e.g. \cite{jones} Corollary 7.1.9.).
 Therefore,
$\kappa:\C(\Z_2\wr H)\to
\cA_H\rtimes H$
extends to a trace (and hence rank) preserving isomorphism between
the von Neumann algebras $\cN(\Z_2\wr H)$ and $\cN(\tau_H)$. 
\begin{propo} \label{elsopropo}
For any countable group $H$,
$$c(\Z_2\wr H)\cong c(\tau_H)\,.$$
\end{propo}
\proof
The rank preserving isomorphism $\kappa:\cN(\Z_2\wr H)\to \cN(\tau_H)$
extends to a rank preserving isomorphism between the rank completions, that is,
the algebras of affiliated operators. It is enough to prove that the rank
closure of $\cA_H\rtimes H$ is $L^\infty_c(B_H,\nu_H)\rtimes H$.

\begin{lemma}
\label{frank}
Let $f\in L^\infty_c(B_H,\nu_H)$. Then $\rk_{\cN(\tau_H)}(f)=\nu_H(\supp (f))$.
\end{lemma}

\proof By definition,
$$
\rk_{\cN(\tau_H)}(f)=1-\lim_{\lambda\to 0} \tr_{\cN(\tau_H)} E_\lambda\,,$$
where $E_\lambda$ is the spectral projection of $f^*f$ corresponding to
$\lambda$. 

$$\tr_{\cN(\tau_H)} E_\lambda=\nu_H(\{x\,\mid\,|f^2(x)|\leq \lambda\})\,.$$
Hence, $\rk_{\cN(\tau_H)}(f)=1-\nu_H(\{x\,\mid f^2(x)=0\})=\nu_H(\supp(f))\,.\quad
\qed$

\noindent
Let $\{m_n\}^\infty_{n=1}\subset \cA_H, m_n\stackrel{\rk}{\to} m\in
L^\infty_c(B_H,\nu_H)$. 
Then $m_n\cdot\gamma \stackrel{\rk}{\to} m\cdot\gamma\,.$
Therefore our proposition follows from the lemma below.
\begin{lemma} \label{density}
$\cA_H$ is dense in $L^\infty_c(B_H,\nu_H)$ with respect to the
rank metric.
\end{lemma}
\proof
By Lemma \ref{frank}, $L^\infty_{fin}(B_H,\nu_H)$ is dense in
$L^\infty_c(B_H,\nu_H)$, where
$L^\infty_{fin}(B_H,\nu_H)$ is the $*$-algebra of functions
taking only finitely many values. Recall that
$V\subset B_H$ is a basic set if $1_V\in\cA_H$. It is well-known
that any measurable set in $B_H$ can be approximated by basic sets,
that is for any $U\subset B_H$, there 
exists a sequence of basic sets $\{V_n\}^\infty_{n=1}$ such 
that
\begin{equation}
\label{e1}
\lim_{n\to\infty} \nu_H(V_n\triangle U)=0\,.
\end{equation}
By (\ref{e1}) and Lemma \ref{frank}
$$\lim_{n\to\infty} \rk_{\cN(\tau_n)}(1_{V_n}-1_U)=0\,.$$

\noindent
Let $f=\sum^l_{m=1} c_m 1_{U_m}$, where $U_m$ are disjoint measurable sets.
Let $\lim_{n\to\infty} \nu_H(V^m_n\triangle U_m)=0$, where
$\{V^m_n\}^\infty_{n=1}$ are basic sets.
Then
$$ \lim_{n\to\infty} \rk_{\cN(\tau_n)} (\sum^l_{m=1} c_m1_{V^m_n}-f)=0\,.$$
Therefore,
$\cA_H$ is dense in $L^\infty_{fin}(B_H,\nu_H)\,.$ \qed

\section{The Odometer Algebra}
\noindent
The Odometer Algebra is constructed via the odometer action using the algebraic
crossed product construction. Let us consider the compact group of $2$-adic
integers $\dya$. Recall that $\dya$ is the completion of the integers
with respect to the dyadic metric
$$d_{(2)}(n,m)=2^{-k}\,,$$
where $k$ is the power of two in the prime factor decomposition of $|m-n|$.
The group $\dya$ can be identified with the compact group of one way infinite
sequences with respect to the binary addition.

\noindent
The Haar-measure $\haar$ on  $\dya$ is defined by $\haar(U^l_n)=1/2^n$,
where $0\leq l \leq 2^n-1$ and $U^l_n$ is the clopen subset of elements in
 $\dya$ having residue $l$ modulo $2^n$.
Let $T$ be the addition map
$x\to x+1$ in $\dya$. The map $T$ defines an action 
$\rho:\Z\curvearrowright (\dya,\haar)$
The dynamical system $(T, \dya,\haar)$
is called the odometer action. As in Section \ref{bernoulli}, we consider
the $*$-subalgebra of function $\cA_M$ in $L^\infty(\dya,\haar)$ that
depend only on finitely many coordinates of $\dya$. We consider
a basis for $\cA_M. $ For $n\geq 0$ and $0\leq l \leq 2^n-1$ let
$$F^l_n(x)=exp\left( \frac{2\pi i x(mod\, 2^n)}{2^n}l \right)\,.$$
Notice that $F^{2l}_{n+1}=F^l_{n}$. Then the functions 
$\{F^l_n\}_{n,l\mid (l,n)=1}$ form the Pr\"ufer $2$-group 
$$\pru=\Z_1\subset\Z_2\subset\Z_4\subset \Z_8\subset\dots$$

\noindent
with respect to the pointwise multiplication. The discrete group $\pru$
is the Pontrjagin dual of the compact Abelian group $\dya$.
The element $F^1_n$ is the generator of the cyclic subgroup $\Z_{2^n}$.
Note that
$$\int_{\dya} F^l_n \,d\haar=0$$
except if $l=0,n=0$, when $F^l_n\equiv 1$. Observe that if $k\in\Z$ then
\begin{equation} \label{bobo1}
\rho(k) F^l_n=F^{l+k(mod\, 2^n)}_n \end{equation}
since $F^l_n(x-k)=F^{l+k(mod\, 2^n)}_n(x)\,.$
Hence we have the following lemma.
\begin{lemma}
The elements of $\cA_M\rtimes \Z$ can be uniquely written as finite sums
in the form 
$$\sum_k \sum_{n\geq 0} \sum_{l\mid (l,n)=1} c_{n,l,k} F^l_n\cdot k\,,$$
where
$k\cdot F^l_n=F^{l+k(mod\, 2^n)}_n$ and $F^0_0=1$.
\end{lemma}

\section{Periodic operators} \label{periodic}
\begin{defin}
A function $\Z\times \Z\to\C$ is a periodic operator if there
exists some $n\geq 1$ such that
\begin{itemize}
\item $A(x,y)=0$, if $|x-y|>2^n$
\item $A(x,y)=A(x+2^n,y+2^n)$.
\end{itemize}
\end{defin}
Observe that the periodic operators form a $*$-algebra, where
\begin{itemize}
\item $(A+B)(x,y)=A(x,y)+B(x,y)$
\item $AB(x,y)=\sum_{z\in \Z} A(x,z)B(z,y)$
\item $A^*(x,y)=A(y,x)$
\end{itemize}
\begin{propo}
The algebra of periodic operators $\cP$ is $*$-isomorphic to a 
dense subalgebra of
$\cM$.
\end{propo}
\proof
We call $A\in\cP$ an element of type-$n$ if
\begin{itemize}
\item $A(x,y)=A(x+2^n,y+2^n)$
\item $A(x,y)=0$ if $0\leq x\leq 2^n-1$, $y>2^n-1$
\item $A(x,y)=0$ if $0\leq x\leq 2^n-1$, $y<0$.
\end{itemize}
Clearly, the elements of type-$n$ form an algebra $\cP_n$ isomorphic to
$\Mat_{2^n\times 2^n}(\C)$ and $\cP_n\to \cP_{n+1}$ is the diagonal embedding.
Hence, we can identify the algebra of finite type elements $\cP_f=
\cup^\infty_{n=1}\cP_n$ with $\varinjlim \Mat_{2^n\times 2^n}(\C)$.

\noindent
For $A\in\cP$, if $n\geq 1$ is large enough,
 let $A_n\in\cP_n$ be defined the following way.
\begin{itemize}
\item $A_n(x,y)=A(x,y)$ if $2^nl\leq x,y\leq 2^nl+2^n-1$ for some $l\in\Z$.
\item Otherwise, $A(x,y)=0$.
\end{itemize}
\begin{lemma} \label{Cauchy}
\begin{description}
\item[(i)] $\{A_n\}^\infty_{n=1}$ is a Cauchy-sequence in $\cM$.
\item[(ii)] $(A+B)_n=A_n+B_n$.
\item[(iii)] $\rk_{\cM} (A_n^*-(A^*)_n)=0$.
\item[(iv)] $\rk_{\cM} ((AB_n)-A_nB_n)=0,.$
\item[(v)] $\lim_{n\to\infty} A_n=0$ if and only if $A=0$.
\end{description}
\end{lemma}

\proof
First observe that for any $Q\in\cP_n$
$$\rk_{\cM} (Q)\leq \frac{|\{ 0\leq x\leq 2^n-1\,\mid \mbox{$\exists\,\,
0\leq y\leq 2^n-1$ such that $A_n(x,y)\neq 0$.}\}|}{2^n}$$

\noindent
Suppose that $A(x,y)=A(x+2^k,y+2^k)$ and $k<n<m$.
Then
$$\left|  \{ 0\leq x\leq 2^n-1\,\mid \, A_n(x,y)\neq A_m(x,y)
\mbox{\,for some $0\leq y\leq 2^n-1$} \} \right| \leq 2^k 2^{m-n}\,.$$

\noindent
Hence by the previous observation, $\{A_n\}^\infty_{n=1}$ is a Cauchy-sequence.
Note that (iii) and (iv) can be proved similarly, 
the proof of (ii) is straightforward.
In order to prove (v) let us suppose
that $A(x,y)=0$ whenever $|x-y|\geq 2^k$.
Let $n>k$ and $0\leq y\leq 2^k-1$ such that $A(x,y)\neq 0$ for some
$-2^k\leq x \leq 2^k-1$.
Therefore $\rk_{\cM} A_n\geq \frac{2^{n-k}-1}{2^n}$.
Thus (v) follows. \qed

\noindent
Let us define $\phi:\cP\to\cM$ by $\phi(A)=\lim_{n\to\infty} A_n$.
By the previous lemma, $\phi$ is an injective $*$-homomorphism. \qed

\begin{defin}
A periodic operator $A$ is diagonal if $A(x,y)=0$, whenever $x\neq y$.
The diagonal operators form the Abelian $*$-algebra $\cD\subset\cP$.
\end{defin}
\begin{lemma}
We have the isomorphism $\cD\cong \C(\pru)$, where $\pru$ is 
the Pr\"ufer $2$-group.
\end{lemma}
\proof
For $n\geq 1$ and $0\leq l \leq 2^n-1$ let $E^l_n\in\cD$ be defined by
$$E^l_n(x,x):=
exp\left( \frac{2\pi i x(mod\, 2^n)}{2^n}l \right)\,.$$
It is easy to see that $E^{2l}_{n+1}=E^l_n$ and
 the multiplicative group generated by $E^1_n$ is isomorphic
to $\Z_{2^n}$. Observe that the set $\{E^l_n\}_{n, l, (l,n)=1}$ form a 
basis in the space of $n$-type diagonal
operators. Therefore, $\cD\cong \cup^\infty_{n=1}\C(Z_{2^n})=\C(\pru)$. \qed

\noindent
Let $J\in\cP$ be the following element.
\begin{itemize}
\item $J(x,y)=1$, if $y=x+1$.
\item Otherwise, $J(x,y)=0$.
\end{itemize}

\noindent
Then 
\begin{equation} \label{bobo2}
J\cdot E^l_n=E^{l+1 (mod\, 2^n)}_n\,. \end{equation}
 Also, any periodic operator
$A$ can be written in a unique way as a finite sum
$$\sum_{k\in\Z} D_k \cdot J^k\,,$$
where $D_k$ is a diagonal operator in the form
$$D_k=\sum^\infty_{n=0} \sum_{l\mid (l,n)=1} c_{l,n,k} E^l_n\,.$$
Thus, by (\ref{bobo1}) and (\ref{bobo2}), we have the following corollary.

\begin{corol} \label{c51} The map $\psi:\cP\to\cA_M\rtimes \Z$
defined by
$$\psi(\sum_k\sum_{n\geq 0}\sum_{l\mid (l,n)=1} c_{l,n,k} E^l_n\cdot k)=
\sum_k\sum_{n\geq 0}\sum_{l\mid (l,n)=1} c_{l,n,k} F^l_n\cdot k $$
is a $*$-isomorphism of algebras. \end{corol}

\section{L\"uck's Approximation Theorem revisited}
The goal of this section is to prove the following proposition.
\begin{propo}
\label{masodikpropo}
We have $c(\rho)\cong \cM$ where $\rho$ is the odometer action.
\end{propo}
\proof
Let us define the linear map $t:\cP\to \C$ by
$$t(A):= \frac{\sum_{i=0}^{2^n-1} A(i,i)}{2^n}\,,$$
where $A\in\cP$ and $A(x+2^n,y+2^n)$ for all $x,y\in \Z$.
\begin{lemma}
\label{traceformula}
$\Tr_{\cN(\rho)}(\psi(A))=t(A)\,,$ where $\psi$ is the $*$-isomorphism
of Corollary \ref{c51}.
\end{lemma}
\proof
Recall that $\Tr_{\cN(\rho)}(F^l_n)=0\,$, except, when $l=0,n=0, F^l_n=1.$
If $n\neq 0$ and $l\neq 0$, then $t(E^l_n)$ is the
sum of all $k$-th roots of unity for a certain $k$, hence $t(E^l_n)=0$. Also,
$t(1)=1$. Thus, the lemma follows. \qed

\vskip 0.2in
\noindent
It is enough to prove that
\begin{equation} \label{egyenlet}
\rk_{\cM}(A)=\rk_{\cN(\rho)}(\psi(A))
\end{equation}

\noindent
Indeed by (\ref{egyenlet}), $\psi$ is a rank-preserving
$*$-isomorphism between $\cP$ and $\cA_M\rtimes \Z$. Hence
the isomorphism $\psi$ extends
to a metric isomorphism
$$\hat{\psi}:\overline{\cP}\to\overline{\cA_M\rtimes \Z}\,,$$
where $\overline{\cP}$ is the closure of $\cP$ in $\cM$ and 
$\overline{\cA_M\rtimes \Z}$ is the closure of
$\cA_M\rtimes \Z$ in $U(\cN(\rho))\,.$
Since $\cP$ is dense in $\cM$, $\overline{\cP}\cong \cM$. Also,
$\overline{\cA_M\rtimes \Z}$ is
a $*$-subalgebra of $U(\cN(\rho))$, since
the $*$-ring operations are continuous with respect to the rank metric.
Therefore $\overline{\cA_M\rtimes \Z}$ is a continuous algebra isomorphic
to $\cM$. Observe that the rank closure $\overline{\cA_M\rtimes \Z}$ 
is isomorphic to the rank closure of $L^\infty_c(\dya,\haar)\rtimes \Z$ by the
argument of Lemma \ref{density}. Therefore,
$c(\rho)\cong\cM$. Thus from now on, our only goal is to prove (\ref{egyenlet}).

\begin{lemma}
\label{elsolemma}
Let $A\in\cP$ and $A_n\in \Mat_{2^n\times 2^n}(\C)$
as in Section \ref{periodic}. Then the matrices $\{A_n\}^\infty_{n=1}$ have
uniformly bounded norms.
\end{lemma}
\proof
Let $M,N$ be chosen in such a way that
\begin{itemize}
\item 
$|A_n(x,y)|\leq M$ for any $x,y\in\Z, n\geq 1$.
\item
$|A_n(x,y)|=0$ if $|x-y|\geq \frac{N}{2}\,.$
\end{itemize}

\noindent
Now let $v=(v(1),v(2),\dots,v(2^n))\in \C^{2^n}$, $\|v\|^2=1\,.$
Then
$$\|A_nv\|^2=\sum^{2^n}_{x=1} |\sum_{y\,\mid |x-y|<N/2} A_n(x,y)v(y)|^2\leq
M^2 \sum^{2^n}_{x=1} |\sum_{y\,\mid |x-y|<N/2} v(y)|^2\leq $$
$$\leq M^2 N \sum^{2^n}_{x=1}\sum_{y\,\mid |x-y|<N/2} |v(y)|^2\leq
M^2  \sum^{2^n}_{y=1} N |v(y)|^2=M^2N^2\,.$$

\noindent
Therefore, for any $n\geq 1$, $\|A_n\|\leq MN\,.$ \qed

\begin{lemma}
\label{masodiklemma}
Let $A\in\cP$. Then for any $k\geq 1$
$$\lim_{k\to\infty} t((A_n^* A_n)^k)=t((A^*A)^k)= 
\Tr_{\cN(\rho)}(\psi(A^*A)^k))\,.$$
\end{lemma}
\proof
Let $m\geq 1, l\geq 1, q\geq 1$ be integers such that
\begin{itemize}
\item
$A(x,y)=A(x+2^m,y+2^m)$ for any $x,y\in\Z$.
\item $A(x,y)=0$, if $|x-y|\geq l$.
\item $|(A^*A)^k(x,x)|\leq q$ and $|(A^*_nA_n)^k(x,x)|\leq q$ for any $x\in\Z$.
\end{itemize}
By definition,
$$ t((A_n^* A_n)^k)=\frac{\sum^{2^n}_{x=1} (A_n^*A_n)^k(x,x)}{2^n}$$
$$ t((A^* A)^k)=\frac{\sum^{2^n}_{x=1} (A^*A)^k(x,x)}{2^n}\,.$$
Observe that if
$2lk<x,2^n-2lk$, then
$$(A^*A)^k(x,x)=(A_n^*A_n)(x,x)\,.$$
Hence,
$$|t((A^* A)^k)-t((A_n^* A_n)^k)|\leq \frac{ 4klq}{2^n}\,.$$
\noindent
Thus our lemma follows. \qed

\noindent
Now, we follow the idea of L\"uck \cite{lueck}.
Let $\mu$ be the spectral measure of $\psi(A)\in\cN(\rho)$. That is
$$\Tr_{\cN(\rho)}f(A^*A)=\int_0^K f(x)\,d\mu(x)\,,$$
for all $f\in C[0,K]$, where $K>0$ is chosen in such a way that
$\mbox{Spec}\,\psi(A^*A)\subset [0,K]$ and
$\|A_n^*A_n\|\leq K$ for all $n\geq 1$.
Also, let $\mu_n$ be the spectral measure of $A_n^*A_n$, that is,
$$t(f(A^*_nA_n))=\int_0^K f(x)\,d\mu_n(x)\,,$$
or all $f\in C[0,K]$. As in \cite{lueck}, we can see
that the measures $\{\mu_n\}^\infty_{n=1}$ converge weakly to $\mu$.
Indeed by Lemma \ref{masodiklemma},
$$\lim_{n\to\infty} t(P(A_n^*A_n))=\Tr_{\cN(\rho)} P(A^*A)\,$$
for any real polynomial $P$, therefore
$$\lim_{n\to\infty} t(f(A_n^*A_n))=\Tr_{\cN(\rho)} f(A^*A)\,$$
for all $f\in C[0,K]$.

\noindent
Since $\rk_\cM(A_n)=\rk_\cM(A_n^*A_n)$ and
$\rk_{\cN(\rho)}(\psi(A))=\rk_{\cN(\rho)})\psi(A^*A))$,
in order to prove (\ref{egyenlet}) it is enough to see that
$$\lim_{n\to\infty} \rk_{\cM}(A_n^*A_n)= \rk_{\cN(\rho)}(\psi(A^*A))\,.$$
Observe that
$ \rk_{\cM}(A_n^*A_n)=1-\mu_n(0)$ and
$$\rk_{\cN(\rho)}(\psi(A^*A))=1-\lim_{\lambda\to 0} \Tr_{\cN(\rho)} E_\lambda=
\mu(0)\,.$$
Hence, our proposition follows from the lemma below (an analogue of L\"uck's
Approximation Theorem).
\begin{lemma}
$\lim_{n\to\infty} \mu_n(0)=\mu(0)\,.$
\end{lemma}
\proof
Let $F_n(\lambda)=\int_0^\lambda \,\mu_n(t)\,dt$ and
$F(\lambda)=\int_0^\lambda \,\mu(t)\,dt$ be the
distribution functions of our spectral measures.
Since $\{\mu_n\}^\infty_{n=1}$ weakly converges to the measure $\mu$, it
is enough to show that $\{F_n\}^\infty_{n=1}$ converges uniformly.
Let $n\leq m$ and $D^n_m:\Mat_{2^n\times 2^n}(\C)\to \Mat_{2^m\times 2^m}(\C)$
be the diagonal operator.
Let $\e>0$. By Lemma \ref{Cauchy}, if $n,m$ are large enough, 
$$\mbox{Rank}(D^n_m(A_n)-A_m)\leq\e 2^m\,.$$
Hence, by Lemma 3.5 \cite{elekrank},
$$\|F_n-F_m\|_\infty\leq \e\,.$$
Therefore, $\{F_n\}^\infty_{n=1}$ converges uniformly. \qed

\section{Orbit Equivalence}
First let us recall the notion of orbit equivalence.
Let $\tau_1:\Gamma_1\curvearrowright(X,\mu)$ resp.
$\tau_2:\Gamma_2\curvearrowright(Y,\nu)$
be essentially free probability measure preserving actions of
the countably infinite groups $\Gamma_1$ resp. $\Gamma_2$.
The two actions are called orbit equivalent if there exists
a measure preserving bijection
$\Psi:(X,\mu)\to (Y,\nu)$ such that
for almost all $x\in X$ and $\gamma\in \Gamma_1$
there exists $\gamma_x\in\Gamma_2$ such that
$$\tau_2(\gamma_x)(\Psi(x))=\Psi(\tau_1(\gamma)(x))\,.$$
Feldman and Moore \cite{feldman} proved
that if $\tau_1$ and $\tau_2$ are orbit equivalent
then $\cN(\tau_1)\cong \cN(\tau_2)\,.$ The goal
of this section is to prove the following proposition.

\begin{propo}
\label{orbit}
If $\tau_1$ and $\tau_2$ are orbit equivalent actions,
then $c(\tau_1)\cong c(\tau_2)$.
\end{propo}

\noindent
Our Theorem \ref{t1} follows from the proposition. Indeed, by Proposition
\ref{elsopropo} and Proposition \ref{masodikpropo}
$$\cM\cong c(\rho)\quad \mbox{and} \quad c(\Z_2\wr H)\cong c(\tau_H)\,.$$
By the famous theorem of Ornstein and Weiss \cite{ornstein}, the odometer
action and the Bernoulli shift action of a countably infinite amenable
group are orbit \\ equivalent. Hence $\cM\cong c(\Z_2\wr H)\,.$\quad\qed
\proof
We build the proof of our proposition on the original proof of Feldman
and Moore.
Let $\gamma\in\Gamma_1$, $\delta\in\Gamma_2$. Let
$$M(\delta,\gamma)=\{y\in Y\,\,\mid \tau_2(\delta)(y)=
\Psi(\tau_1(\gamma)\Psi^{-1}(y))\}$$
$$N(\gamma,\delta)=\{x\in X\,\,\mid \tau_1(\gamma)(x)=
\Psi^{-1}(\tau_2(\delta)\Psi(x))\}.$$

\noindent
Observe that $\Psi(N(\delta,\gamma))=M(\gamma,\delta)\,.$
following Feldman and Moore (\cite{feldman}, Proposition 2.1)
for any $\gamma\in\Gamma_1$, $\delta\in\Gamma_2$
$$\kappa(\gamma)=\sum_{h\in\Gamma_2} h\cdot 1_{M(h,\gamma)}$$
and
$$\lambda(\delta)=\sum_{g\in\Gamma_1} g\cdot 1_{N(g,\delta)}$$
are well-defined. That is,
$\sum^k_{n= 1} h_n\cdot 1_{M(h_n,\gamma)}$ converges weakly 
to $\kappa(\gamma)\in\cN(\tau_2)$ as $k\to\infty$ and
$\sum^k_{n= 1} g_n\cdot 1_{N(g_n,\delta)}$ converges weakly 
to $\lambda(\delta)\in\cN(\tau_1)$ as $k\to\infty$, where
$\{\gamma_n\}^\infty_{n=1}$ resp. $\{\delta_n\}^\infty_{n=1}$
are enumerations of the elements of $\Gamma_1$ resp. $\Gamma_2$.

\noindent
Furthermore, one can extend $\kappa$ resp. $\lambda$ to maps
$$\kappa':L^\infty((X,\mu)\rtimes\Gamma_1)\to \cN(\tau_2)$$ resp.
$$\lambda':L^\infty((Y,\nu)\rtimes\Gamma_2)\to \cN(\tau_1)$$
by
$$\kappa'(\sum_{\gamma\in\Gamma_1}a_\gamma\cdot\gamma)=
\sum_{\gamma\in\Gamma_1}(a_\gamma\circ\Psi^{-1})\cdot \kappa(\gamma)
=\sum_{\gamma\in\Gamma_1}(a_\gamma\circ\Psi^{-1})\cdot
\sum^\infty_{n=1}h_n\cdot 1_{M(h_n,\gamma)}$$
and
$$\lambda'(\sum_{\delta\in\Gamma_2}b_\delta\cdot\delta)=
\sum_{\delta\in\Gamma_2}(b_\delta\circ\Psi)\cdot \lambda(\delta)
=\sum_{\delta\in\Gamma_2}(b_\delta\circ\Psi)\cdot
\sum^\infty_{n=1}g_n\cdot 1_{N(g_n,\delta)}\,.$$

\noindent
The maps $\kappa'$ resp. $\lambda'$ are injective
trace-preserving $*$-homomorphisms with weakly dense ranges. Hence they
extend
to isomorphisms of von Neumann algebras
$$\hat{\kappa}:\cN(\tau_1)\to \cN(\tau_2), 
\hat{\lambda}:\cN(\tau_2)\to \cN(\tau_1)\,,$$
where
$\hat{\kappa}$ and $\hat{\lambda}$ are, in fact, the inverses of each other.
\begin{lemma}
\begin{equation}
\label{rank1}
\lim_{k\to\infty} \rk_{\cN(\tau_2)}\left(
\sum_{\gamma\in\Gamma_1}(a_\gamma\circ\Psi^{-1})\cdot
\sum^k_{n=1}h_n\cdot 1_{M(h_n,\gamma)}-
\hat{\kappa}(\sum_{\gamma\in\Gamma_1}a_\gamma\cdot\gamma)
\right)=0\,.
\end{equation}
\begin{equation}
\label{rank2}
\lim_{k\to\infty} \rk_{\cN(\tau_1)}\left(
\sum_{\delta\in\Gamma_2}(b_\delta\circ\Psi)\cdot
\sum^k_{n=1}g_n\cdot 1_{N(g_n,\delta)}-
\hat{\lambda}(\sum_{\delta\in\Gamma_2} b_\delta\cdot\delta)\right)=0\,.
\end{equation}
\end{lemma}

\proof
By definition, the disjoint union $\cup_{n=1}^\infty M(h_n,\gamma)$ equals
to $Y$ (modulo a set of measure zero).
We need to show that if
$\{\sum^k_{n=1} T_n\cdot  1_{M(h_n,\gamma)}\}^\infty_{k=1}$ weakly converges
to an element $S\in \cN(\tau_2)$, then
$\{\sum^k_{n=1} T_n\cdot  1_{M(h_n,\gamma)}\}^\infty_{k=1}$ converges
to $S$ in the rank metric as well, where 
$T_n\in L^\infty_c(Y,\nu)\rtimes\Gamma_2$.
Let $P_k=\sum^k_{n=1}  1_{M(h_n,\gamma)}\in l^2(\Gamma,L^2(Y,\nu))$.
We denote by $\hat{P}_k$ the element
$\sum^k_{n=1}  1_{M(h_n,\gamma)}$ in $L^\infty_c(Y,\nu)\rtimes\Gamma_2$.
By definition, if $L(A)(P_k)=0$ then $A\hat{P}_k=0$.
Now, by weak convergence,
$$L(S)(P_k)=\lim_{l\to\infty} \sum^l_{n=1} T_n\cdot 1_{M(h_n,\gamma)} (P_k)\,.$$
That is,
$$L(S- \sum^k_{n=1} T_n\cdot 1_{M(h_n,\gamma)})(P_k)=0\,.$$
Therefore,
$$(S- \sum^k_{n=1} T_n\cdot 1_{M(h_n,\gamma)})\hat{P}_k=0\,.$$
Thus,
$$(S- \sum^k_{n=1} T_n\cdot 1_{M(h_n,\gamma)})=
(S- \sum^k_{n=1} T_n\cdot 1_{M(h_n,\gamma)})(1-\hat{P}_k)\,.$$
By Lemma \ref{frank},
$\rk_{\cN(\tau_2)}(1-\hat{P}_k)=1-\sum^k_{n=1} \nu(M(h_n,\gamma))$, hence
$$\lim_{k\to\infty}\rk_{\cN(\tau_2)}
(S- \sum^k_{n=1} T_n\cdot 1_{M(h_n,\gamma)})=0\,.\quad\qed$$

\noindent
Now let us turn back to the proof of our proposition.
By (\ref{rank1}), $\hat{\kappa}$ maps the
algebra $L^\infty_c(X,\mu)\rtimes \Gamma_1$ into
the rank closure
of  $L^\infty_c(Y,\nu)\rtimes \Gamma_2$. Since $\hat{\kappa}$ preserves
the rank, $\hat{\kappa}$ maps the rank closure
of $L^\infty_c(X,\mu)\rtimes \Gamma_1$ into the  rank closure
of  $L^\infty_c(Y,\nu)\rtimes \Gamma_2$.
Similarly, $\hat{\lambda}$ maps 
 the rank closure
of $L^\infty_c(Y,\nu)\rtimes \Gamma_2$ into the  rank closure
of  $L^\infty_c(X,\mu)\rtimes \Gamma_1$. That is,
$\hat{\kappa}$ provides an isomorphism
between the rank closures
of $L^\infty_c(X,\mu)\rtimes \Gamma_1$ and $L^\infty_c(Y,\nu)\rtimes \Gamma_2$.
Therefore, the smallest continuous ring containing 
$L^\infty_c(X,\mu)\rtimes \Gamma_1$ in $U(\cN(\tau_1))$ is mapped to 
the smallest continuous
ring containing $L^\infty_c(Y,\nu)\rtimes \Gamma_2$ in $U(\cN(\tau_2))$ . \qed

\noindent
G\'abor Elek \\
Lancaster University \\
Department of Mathematics and Statistics \\
g.elek@lancaster.ac.uk

\end{document}